\documentclass[a4paper,12pt]{article}

\usepackage{amsmath}
\usepackage{amsthm}
\usepackage{amsfonts}

\parskip\smallskipamount
\newtheorem*{theorem*}{Theorem}
\newtheorem{theorem}{Theorem}

\newtheorem{remark}{Remark}

\newcommand{\E}{\mathbb{E}}

\newcommand{\Z}{\mathbb{Z}}

\newcommand{\Na}{\mathbb{N}}

\title{A note on spatial monotonicity for \\ one-dimensional spin systems  \vspace{-3mm}}
\date{}
\author{Tzioufas, A.\footnote{On leave from the University of Buenos Aires}} 
\begin{document}
\maketitle
\vspace{-10mm}

\begin{abstract} 
We show that attractive, translation invariant, one-dimensional spin systems started from the unit step function possess the following basic property. At any time, the entire configuration from a point onward is stochastically decreasing with respect to distance from this point to the origin. The proof relies on a stochastic domination argument which exploits the interaction assumptions in a simple manner. 
\end{abstract}

%The probability that a site is included in the process decreases monotonically with respect to its distance from the origin.

%\subsection{Attractive spin systems}

{\bf\noindent Introduction. } We consider \textit{one-dimensional spin systems with attractive, translation invariant, and finite-range interaction}. These are Feller processes on the compact state space $\{0,1\}^{\Z}$ endowed with the product topology. These processes constitute an extensive, well-developed field of study, containing various important special cases, for instance, contact and voter processes, or stochastic Ising models; for overview material cf.\ with Liggett \cite{L10}, Chpt.\ 4; Liggett \cite{L85}, Chpt.\ III; and Griffeath \cite{G79}. 

From the numerous interpretations of the coordinate values, which depend on the application in mind, we adopt that of the evolution of a population of particles whereby the value "1" corresponds to {\it occupied } 
site, whereas that of "0" corresponds to {\it empty}. Thus, the process is thought of as modeling the one-dimensional spatial growth associated with the set of occupied sites. Spin systems evolve according to transitions that involve change of only one coordinate 
with the restriction that birth attempts to already occupied sites being suppressed, i.e.\ sites are able to accommodate no more than one particle. Their definition is facilitated by the so-called flip rates, denoted by $c(x,\eta)$ here and below, which is a translation invariant, non-negative, bounded function on $\Z \times \{0,1\}^{\Z}$ that, informally speaking, corresponds to the rate at which the coordinate, or spin, $\eta(x)$ flips from 0 to 1 or from 1 to 0 when the state of the process is $\eta$. The assumption of finite range refers to that $c(x,\eta)$ may depend on $\eta$ only through the spin at sites in an arbitrary but finite neighborhood of $x$, and is known to ascertain the uniqueness and well-definedness of the process (cf., for instance, Theorem 3.4 in \cite{L10}). Attractiveness (cf.\ (\ref{attr}) below) ensures that the process is monotone in the sense that stochastic order is preserved by the evolution, and is a fundamental, amply exploited  assumption (cf.\ Chpt.\ III.2 in \cite{L85}) that is patently natural (cf.\ (\ref{nsattr}) below) and weaker than another property called additivity (cf.\ with Corollary 1.3,  Chpt.\ 2, in Griffeath \cite{G79}).  For an extensive list of classes of processes which fulfill these standard assumptions, in addition to those mentioned above, we refer to the example paragraphs in Chpt.\ III, \cite{L85}.

% in the standard particle births and deaths interpretation of the process
%between this type of coupling and coalescing duality of additive, finite-range spin systems
% or the set of wet sites on the diagonal in percolation on the so-called "north-east" lattice, see the figure in p.\ 367 in \cite{G10}
%\vspace{2mm}

Connectivity functions are principal objects in percolation theory that concern the probability of connecting the origin by open paths to collections of lattice points as a function of their relative position to the origin (cf.\ with Grimmett \cite{G10}). These functions possess natural counterparts for spin systems which admit a graphical, or percolation substructure, representation (cf.\ with for instance Liggett \cite{L85}, Chpt.\ III.6 for details about this important link). Basic properties of these functions are therefore of inherent interest, as these are fundamental objects in the understanding and analysis of these processes, in addition to that it is natural to expect for this same reason their involvement in other arguments and in applications. 
In this regard, soon after Bezuidenhout and Grimmett \cite{BG90},  Gray \cite{G91} raises three questions regarding  one of the central spin-systems, the \textit{basic contact process}\footnote{that is, a one-dimensional (finite) nearest neighbors and symmetric interaction particle system, which in discrete-time corresponds to standard two-dimensional oriented percolation (cf. Durrett \cite{D84}).}, viewed there as so fundamental that the lack of their treatment at the time is perceived there as "somewhat embarrassing".  The first of these questions comprises the basic and intuitively plausible property that, when the process is started from only the origin occupied, at any time, the probability of a site being occupied decreases monotonically with respect to its distance from the origin.  To see the connection, bearing in mind that occupancy corresponds to existence of an open path in the graphical representation, note that this property may also be rendered equivalently as a monotonicity property of the associated two-point connectivity functions. We refer to Theorem 3 in \cite{G91} for precise formulations and a proof outline and, further, to Andjel and Gray \cite{AG14} and Andjel and Sued \cite{AS08}, where detailed proofs may be found, and note that these proofs rely crucially on path intersection properties available only under the nearest neighbor interaction assumption. \footnote{The second and third questions raised there regard a non-intuitive correlation inequality for the events of two sites with a certain spatial ordering being occupied, and another property, conjectured in Liggett \cite{L86}, regarding the survival probabilities of processes started from configurations with an interval gap (see Theorems 4 and 5 there; see also \cite{AG14} for a proof of the the former result).}

Another look is taken here in regard to extensions of the property comprising the first of the aforementioned questions over attractive, finite-range, one-dimensional spin systems. We show that when started from the unit step function, i.e.\ sites are occupied if and only if are non-positive, these processes fulfill collectively the following stronger property.  At any time, the joint distribution of the \textit{entire configuration} of the process from a point onwards is stochastically decreasing with respect to this points' distance from the origin. The result is formally presented in the next section, and its proof may be found in the immediately following one. 
\vspace{2mm}

{\bf\noindent Statement of Result. }Let $c(x,\eta)$ be a translation invariant, non-negative and bounded function on $\Z \times \Xi$, $\Xi :=\{0,1\}^{\Z}$, 
which is of \textit{finite-range}, that is, depends on $\eta$ only through its value at $x$ and a finite number of its coordinates. 
The one-dimensional spin system $\xi_{t}$ with flip rates $c(x,\eta)$ is the Feller process on $\Xi$ with generator $\mathcal{L}$ given by  
\begin{equation}\label{El}
\mathcal{L} f (\eta) = \sum_{x} c(x, \eta) ( f(\eta_{x})- f(\eta)), 
\end{equation}
for $f$ that depend on finitely many coordinates, and where $\eta_{x}$ is such that $\eta_{x}(x) = 1-\eta(x)$ and $\eta_{x}(y) = \eta(y)$ for all $y \not= x$.  For a formal construction and regarding the existence and uniqueness of the process, ascertained under our assumptions, see for instance Theorem B.3 in \cite{L99}. We will also adapt the assumption of \textit{attractiveness} which provides that: 
\begin{equation}\label{attr}
\mu_{1} \leq \mu_{2} \Rightarrow \mu_{1}S(t) \leq \mu_{2}S(t)
\end{equation}
for all  $t\geq0$, where $S(t)$ denotes the associated semigroup of the process, and the usual partial order on the set of probability measures on $\Xi $ is in force (that is, if $\mu$ and $\nu$ are probability measures on $\Xi $ we write $\mu \leq \nu$ whenever $\int f d\mu \leq \int f d\nu$, for all $f$ bounded and increasing w.r.t.\ the coordinatewise partial order on $\Xi$). The following is a simple necessary and sufficient condition for attractiveness. For any $\eta, \zeta$ such that $\eta \leq \zeta$, we have that 
\begin{align}\label{nsattr} 
\mbox{if } \eta(x) = \zeta(x) = 0, & \mbox{ then } c(x,\eta) \leq c(x, \zeta), \nonumber \\ 
\mbox{if }  \eta(x) = \zeta(x) = 1, & \mbox{ then } c(x,\eta) \geq c(x, \zeta), 
\end{align}
for a proof see Theorem 2.2, p.\ 134 \cite{L85}.  In the particle population interpretation adopted, note that this condition comprises of birth rates being increasing functions of the configuration and death rates being decreasing. 

%, thinking of the flips from $0$ to $1$ as births of particles and those from $1$ to $0$ as deaths,

%the space of continuous functions on $\Xi $ with the uniform norm
% Accordingly, a spin system is said to be attractive, finite range when its flip rates satisfy the corresponding properties.  of $\eta$,

%,

%Introduce oredering.

\noindent\textbf{Definition.} Let $\Na ^{+}$ denote the non-negative integers. We say that a collection of probability measures $(\nu_{n}: n\geq0)$ on $\Xi^{+} := \{0,1\}^{\Na^{+}}$ is \textit{decreasing} whenever 
\[
\nu_{n} \geq \nu_{n+1}, \hspace{3mm} \mbox{ for all }   n \geq0,
\] 
according to the usual partial order of probability measures on $\Xi^{+}$ (see that elaborated after (\ref{attr})).

To state the result let also $1_{\cdot}$ denote the indicator function.   

%Let $\eta^{n}_{k,l}$ be the configuration in $\Xi^{+}$ such that $\eta^{n}_{k,l}(x) = 1$ if and only if $-k-1 \leq x \leq 0$.  
%for all large enough $n\geq0$.

%\xi_{t=0}^{(-\infty,0]} =
\begin{theorem}\label{tspin}
Let $\bar{\xi}_{t}$ be any attractive and finite-range one-dimensional spin system with $\bar{\xi}_{t=0} = 1_{(-\infty,0]}$. Let also $\mu_{t}(z)$ denote the joint law of $(\bar{\xi}_{t}(z), \bar{\xi}_{t}(z+1), \dots)$. For all $t$, we have that $(\mu_{t}(z): z\geq0)$ is decreasing.  
\end{theorem}

\begin{remark}\label{rem2}  \normalfont  
By an elementarily limits procedure, Theorem \ref{tspin} extends to $\xi_{t}'$ with $\xi'_{t=0} = 1_{[-N,0]}$, for all sufficiently large, but finite $N$. This may be useful in conjunction with renormalization arguments (cf.\ Bezuidenhout and Grimmett \cite{BG90} and  Bezuidenhout and Gray \cite{BG94}). %B) Works for any periodic initial configuration - the 1111.. have period 0.
\end{remark}
%We now give the proof of Theorem \ref{tspin}. 
\begin{proof}[{\bf Proof of Theorem \ref{tspin}}]
Given configuration $\eta$, let $\eta[z, \infty)$ denote the element of $\Xi^{+}$ such that 
\[
\eta[z, \infty)(x) := \eta(x+z), \hspace{3mm} \mbox{for all } x \geq0.
\] 
We will show that 
\begin{equation}\label{finf}
\bar{\xi}_{t}[z, \infty) \geq_{st.} \bar{\xi}_{t}[z+1, \infty)
\end{equation}
where  $\geq_{st.}$ denotes stochastic domination, defined here as
\[
\E[f\big(\bar{\xi}_{t}[z, \infty)\big)] \geq  \E[f \big(\bar{\xi}_{t}[z+1, \infty) \big)],  
\]
for all bounded, non-decreasing $f: \Na \rightarrow \{0,1\}$.

% monotonicity in the initial configuration, which is a consequence of the coupling provided by the graphical representation,

Let $\tilde{\xi}_{t}$ denote the process with $\tilde{\xi}_{t=0} = 1_{(-\infty,-1]}$. Attractiveness provides that
\begin{equation}\label{eqmon}
\bar{\xi}_{t}[z,\infty) \geq_{st.} \tilde{\xi}_{t}[z,\infty), \mbox{ for all } t. 
\end{equation}
Furthermore, translation invariance ensures that $\bar{\xi}_{t}[z+1,\infty)$ is equal in distribution to $\tilde{\xi}_{t}[z, \infty)$, which combined with $(\ref{eqmon})$ gives $(\ref{finf})$, hence completing the proof. 

\end{proof} 

{\bf Acknowledgements:}  Thanks for related discussions and advice to Pablo Ferrari and Ioannis Kontoyiannis.

\end{document}